\documentclass[12pt]{article}
\usepackage{amsfonts}
\usepackage[latin1]{inputenc}

\usepackage[all]{xy}

\def\squarebox#1{\hbox to #1{\hfill\vbox to #1{\vfill}}}
\newcommand{\qed}{\hspace*{\fill}
\vbox{\hrule\hbox{\vrule\squarebox{.667em}\vrule}\hrule}\smallskip}
\newtheorem{teorema}{Theorem}[section]
\newtheorem{lema}[teorema]{Lemma}
\newtheorem{corolario}[teorema]{Corollary}
\newtheorem{proposicao}[teorema]{Proposition}

\newenvironment{profe}{\noindent {\bf Proof:}}{\hfill $\qed $ \newline}

\begin{document}

\title{Orientability of vector bundles over real flag manifolds\thanks{%
Research supported by FAPESP grant no.\ 07/06896-5 }}
\author{Mauro Patr\~{a}o\thanks{%
Department of Mathematics. Universidade de Bras\'{\i}ília. Bras\'{\i}lia -
DF, Brazil. } \and Luiz A. B. San Martin\thanks{%
Institute of Mathematics. Universidade Estadual de Campinas. Campinas - SP,
Brazil. } \thanks{%
Supported by CNPq grant no.\ 303755/09-1} \and La\'{e}rcio J. dos Santos%
\thanks{%
Institute of Mathematics. Universidade Federal de S\~{a}o Carlos. Sorocaba -
SP, Brazil. } \and Lucas Seco\thanks{%
Department of Mathematics. Universidade de Brasília. Bras\'{\i}lia -
DF, Brazil. }}
\date{}
\maketitle

\begin{abstract}
We investigate the orientability of a class of vector bundles over flag manifolds of real semi-simple Lie groups, which include the tangent bundle and also stable bundles of certain gradient flows. Closed formulas, in terms of roots, are provided.  
\end{abstract}



\section{Introduction}

We investigate the orientability of a class of vector bundles over flag manifolds of real semi-simple Lie groups, the so called (generalized) real flag manifolds.  These include the tangent bundle and also stable bundles of some gradient flows on these manifolds which were considered elsewhere (see Section 3 of Duistermaat-Kolk-Varadarajan \cite{dkv} and Section \ref{secbundlesbundles} of the present article).  We get closed formulas, in terms of roots associated to the real flag manifolds, to decide when they are orientable.  As far as we know,  our results and methods of proof are not known.

The topology of flag manifolds of complex semi-simple Lie groups, and of holomorphic vector bundles over them is, by now, a well-understood classical subject (see, for example, Bernstein-Gel'fand-Gel'fand \cite{bgg} or Bott-Borel-Weil's Theorem \cite{bbw}).  On the other hand, the topology of real flag manifolds is a more delicate subject.  Its mod 2 homology was obtained in the 1980's (see Section 4 of \cite{dkv}) and in the 1990's it was obtained a complete (although algorithmic) description of its integral homology (Kocherlakota \cite{koch}, see also \cite{losm}) and its fundamental groups (Wiggerman \cite{wiggerman}).  It is begining to emerge relations between the cohomology of real flag manifolds and infinite dimensional representation theory of the real semi-simple Lie group (Casian-Stanton \cite{casian1}) and dynamics of integrable systems (Casian-Kodama \cite{casian2, casian3}).  As for the topology of vector bundles over real flag manifolds, we are not aware of any general result in the literature. This article is a contribution in this direction.

The structure of the article is as follows.  In Section \ref{secprelim} we recall some definitions and facts about real semi-simple Lie groups and their flag manifolds. In particular we look at the structure of the connected components certain centralizers that will appear later as isotropy subgroups (Subsection \ref{compon}). Also we recall  the construction of the stable and unstable vector bundles over fixed points of gradient flows (Subsection \ref{secfibr}).  For these stable bundles, and also for the tangent bundle of a real flag manifold, there is a Lie group acting on the vector bundle by linear maps in such a way that the action on the base space is transitive. In both cases, the base space is a homogenous space of a Lie group.

In Section \ref{secvechomog} we derive our method of determining orientability of vector bundles over a homogeneous space of a Lie group, which consists of reducing the
orientability question to a computation of signs of determinants.  Namely the
vector bundle is orientable if and only if each linear map coming from the
representation of the isotropy subgroup on the fiber at the origin has
positive determinant (see Proposition \ref{proporient}).  Using this criterion we get closed formulas, in terms of roots and their multiplicities to decide when one of our vector bundles is orientable (see Theorems \ref{teoorientflags} and \ref{teoorientstable}, below). In particular, we prove that any maximal flag manifold is orientable. A
result already obtained by Kocherlakota \cite{koch} as a consequence of the
computation of the homology groups of the real flag manifolds. 

In Section \ref{secsplit} we make a detailed analysis of the orientability of
the flag manifolds associated to the split real forms of the classical Lie
algebras $A_{l}=\frak{sl}\left( l+1,\mathbb{R}\right) $, $B_{l}=\frak{so}%
\left( l,l+1\right) $, $C_{l}=\frak{sp}\left( l,\mathbb{R}\right) $ and $D_{l}=%
\frak{so}\left( l,l\right) $.

The orientability of the stable and unstable bundles was our original
motivation to write this paper. It comes from the computation of the Conley
indices for flows on flag bundles in \cite{psmsconley}. In this computation
one wishes to apply the Thom isomorphism between homologies of the base
space and the disk bundle associated to a vector bundle. The isomorphism
holds in $\mathbb{Z}$ homology provided the bundle is orientable, asking for
criteria of orientability of such bundles. We develop along this line on
Section \ref{secbundlesbundles}.

\section{Preliminaries}\label{secprelim}

We recall some facts of semi-simple Lie groups and their flag manifolds (see
Duistermat-Kolk-Varadarajan \cite{dkv}, Helgason \cite{he}, Humphreys \cite{hu} Knapp \cite{knp}
and Warner \cite{w}). To set notation let $G$ be a connected noncompact real
semi-simple Lie group with Lie algebra $\frak{g}$. Fix a Cartan involution $%
\theta $ of $\frak{g}$ with Cartan decomposition $\frak{g}=\frak{k}\oplus 
\frak{s}$. The form $\langle X, Y\rangle_{\theta} =-\langle X,\theta
Y\rangle $, where $\langle \cdot ,\cdot \rangle $ is the Cartan-Killing form
of $\frak{g}$, is an inner product. An element $g \in G$ acts in $X \in 
\frak{g}$ by the adjoint representation and this is denoted by $gX$.

Fix a maximal abelian subspace $\frak{a}\subset \frak{s}$ and a Weyl chamber 
$\frak{a}^{+}\subset \frak{a}$. We let $\Pi $ be the set of roots of $\frak{a%
}$, $\Pi ^{+}$ the positive roots corresponding to $\frak{a}^{+}$, $\Sigma $
the set of simple roots in $\Pi ^{+}$ and $\Pi^- = - \Pi^+$ the negative
roots. The Iwasawa decomposition of the Lie algebra $\frak{g}$ reads $\frak{g%
}=\frak{k}\oplus \frak{a}\oplus \frak{n}^{\pm}$ with $\frak{n}%
^{\pm}=\sum_{\alpha \in \Pi ^{\pm}}\frak{g}_{\alpha }$ where $\frak{g}%
_{\alpha }$ is the root space associated to $\alpha $. As to the global
decompositions of the group we write $G=KS$ and $G=KAN^\pm$ with $K=\exp 
\frak{k}$, $S=\exp \frak{s}$, $A=\exp \frak{a}$ and $N^{\pm}=\exp \frak{n}%
^{\pm}$.

The Weyl group $W$ associated to $\frak{a}$ is the finite group generated by
the reflections over the root hyperplanes $\alpha =0$ in $\frak{a}$, $\alpha
\in \Pi $. $W$ acts on $\frak{a}$ by isometries and can be alternatively be
given as $W=M^{*}/M$ where $M^{*}$ and $M$ are the normalizer and the
centralizer of $A$ in $K$, respectively. We write $\frak{m}$ for the Lie
algebra of $M$.

\subsection{Subalgebras defined by simple roots}

Associated to a subset of simple roots $\Theta \subset \Sigma $ there are
several Lie algebras and groups (cf. \cite{w}, Section 1.2.4): We write $%
\frak{g}\left( \Theta \right) $ for the (semi-simple) Lie subalgebra
generated by $\frak{g}_{\alpha }$, $\alpha \in \Theta $, put $\frak{k}%
(\Theta )=\frak{g}(\Theta )\cap \frak{k}$ and $\frak{a}\left( \Theta \right)
=\frak{g}\left( \Theta \right) \cap \frak{a}$. The simple roots of $\frak{g}%
(\Theta )$ are given by $\Theta $, more precisely, by restricting the
functionals of $\Theta $ to $\frak{a}(\Theta )$. Also, the root spaces of $%
\frak{g}(\Theta)$ are given by $\frak{g}_\alpha$, for $\alpha \in \langle
\Theta \rangle$. Let $G\left( \Theta \right) $ and $K(\Theta )$ be the
connected groups with Lie algebra, respectively, $\frak{g}\left( \Theta
\right) $ and $\frak{k}\left( \Theta \right) $. Then $G(\Theta )$ is a
connected semi-simple Lie group.

Let $\frak{a}_{\Theta }=\{H\in \frak{a}:\alpha (H)=0,\,\alpha \in \Theta \}$
be the orthocomplement of $\frak{a}(\Theta )$ in $\frak{a}$ with respect to
the $\langle \cdot, \cdot \rangle_{\theta}$-inner product. We let $K_{\Theta
}$ be the centralizer of $\frak{a}_{\Theta }$ in $K$. It is well known that 
\[
K_{\Theta } = M(K_{\Theta })_0 = MK(\Theta ). 
\]
Let $\frak{n}_{\Theta }^{\pm }=\sum_{\alpha \in \Pi ^{\pm }-\langle \Theta
\rangle }\frak{g}_{\alpha }$ and $N_{\Theta }^{\pm }=\exp (\frak{n}_{\Theta
}^{\pm })$. We have that $K_\Theta$ normalizes $\frak{n}_{\Theta }^{\pm }$
and that $\mathfrak{g}=\frak{n}_{\Theta }^{-}\oplus \frak{p}_{\Theta }$. The
standard parabolic subalgebra of type $\Theta \subset \Sigma$ with respect
to chamber $\frak{a}^+$ is defined by 
\[
\frak{p}_{\Theta }=\frak{n}^{-}\left( \Theta \right) \oplus \frak{m}\oplus 
\frak{a}\oplus \frak{n}^{+}. 
\]
The corresponding standard parabolic subgroup $P_{\Theta }$ is the
normalizer of $\frak{p}_{\Theta }$ in $G$. It has the Iwasawa decomposition $%
P_{\Theta } = K_\Theta A N^+$. The empty set $\Theta =\emptyset $ gives the
minimal parabolic subalgebra $\frak{p} = \frak{m}\oplus \frak{a}\oplus \frak{%
n}^{+}$ whose minimal parabolic subgroup $P=P_\emptyset $ has Iwasawa
decomposition $P=MAN^+$.

Let $d = \mathrm{dim}( \frak{p}_{\Theta } )$ and consider the Grassmanian of 
$d$-dimensional subspaces of $\frak{g}$, where $G$ acts by its adjoint
representation. The flag manifold of type $\Theta $ is the $G$-orbit of the
base point $b_{\Theta }=\frak{p}_{\Theta }$, which we denote by $\Bbb{F}%
_{\Theta }$. This orbit identifies with the homogeneous space $G/P_{\Theta }$%
. Since the adjoint action of $G$ factors trough $\mathrm{Int}(\mathfrak{g})$%
, it follows that the flag manifolds of $G$ depends only on its Lie algebra $%
\frak{g}$. The empty set $\Theta =\emptyset $ gives the maximal flag
manifold $\Bbb{F}=\Bbb{F}_{\emptyset }$ with basepoint $b=b_{\emptyset }$.

\subsection{Subalgebras defined by elements in $\frak{a}$}

The above subalgebras of $\frak{g}$, which are defined by the choice of a
Weyl chamber of $\frak{a}$ and a subset of the associated simple roots, can
be defined alternatively by the choice of an element $H \in \frak{a}$ as
follows. First note that the eigenspaces of $\mathrm{ad}(H)$ in $\frak{g}$
are the weight spaces $\frak{g}_\alpha$. Now define the negative and
positive nilpotent subalgebras of type $H$ given by 
\[
\frak{n}^-_H = \sum\{ \frak{g}_\alpha:\, \alpha(H) < 0 \}, \qquad \frak{n}%
^+_H = \sum\{ \frak{g}_\alpha:\, \alpha(H) > 0 \}, 
\]
and the parabolic subalgebra of type $H$ which is given by 
\[
\frak{p}_H = \sum\{ \frak{g}_\alpha:\, \alpha(H) \geq 0 \}. 
\]
Denote by $N^\pm_H = \exp(\frak{n}^\pm_H)$ and by $P_H$ the normalizer in $G$
of $\frak{p}_H$. Let $d = \mathrm{dim}( \frak{p}_H )$ and consider the
Grassmanian of $d$-dimensional subspaces of $\frak{g}$, where $G$ acts by
its adjoint representation. The flag manifold of type $H$ is the $G$-orbit
of the base point $\frak{p}_H$, which we denote by $\Bbb{F}_H$. This orbit
identifies with the homogeneous space $G/P_H$, where $P_H$ is the normalizer
of $\frak{p}_H$ in $G$.

Now choose a chamber $\frak{a}^+$ of $\frak{a}$ which contains $H$ in its
closure, consider the simple roots $\Sigma$ associated to $\frak{a}^+$ and
consider 
\[
\Theta (H)=\{\alpha \in \Sigma :\alpha (H)=0\}, 
\]
the set of simple roots which annihilate $H$. Since a root $\alpha \in
\Theta(H)$ if, and only if, $\alpha|_{\frak{a}_{\Theta(H)}} = 0$, we have
that 
\[
\frak{n}^\pm_H = \frak{n} ^\pm_{\Theta(H)} \quad\mbox{and}\quad \frak{p}_H = 
\frak{p}_{\Theta(H)}. 
\]
Denoting by $K_H$ the centralizer of $H$ in $K$, we have that $K_{H}=
K_{\Theta (H)}$. So it follows that 
\[
{\mathbb F}_H = {\mathbb F}_{\Theta(H)}, 
\]
and that the isotropy of $G$ in $\frak{p}_H$ is 
\[
P_H = P_{\Theta(H)} = K_{\Theta(H)} A N^+ = K_H A N^+, 
\]
since $K_{\Theta(H)} = K_H$. Denoting by $G(H) = G(\Theta(H))$ and by $K(H)
= K(\Theta(H))$, it is well know that 
\[
K_H = M (K_H)_0 = M K(H). 
\]

We remark that the map 
\begin{equation}  \label{eq:mergulho-em-s}
{\mathbb F}_H \to \mathfrak{s},\qquad k \frak{p}_H \mapsto k H,\quad%
\mbox{where } k \in K,
\end{equation}
gives an embeeding of ${\mathbb F}_H$ in $\mathfrak{s}$ (see Proposition 2.1
of \cite{dkv}). In fact, the isotropy of $K$ at $H$ is $K_H = K_{\Theta(H)}$
which is, by the above comments, the isotropy of $K$ at $\frak{p}_H$.

\subsection{Connected components of $K_H$}\label{compon}

We assume from now on that $G$ is the adjoint group $\mathrm{Int}\left( 
\frak{g}\right) $. There is no loss of generality in this assumption because
the action on the flag manifolds of any locally isomorphic group factors
through $\mathrm{Int}\left( \frak{g}\right) $. The advantage of taking the
adjoint group is that it has a complexification $G_{\Bbb{C}}=\mathrm{Aut}%
_{0}\left( \frak{g}_{\Bbb{C}}\right) $ with Lie algebra $\frak{g}_{\Bbb{C}}$
in such a way that $G$ is the connected subgroup of $G_{\Bbb{C}}$ with Lie
algebra $\frak{g}$.

For a root $\alpha $, let ${\alpha }^{\vee }=2{\alpha }/\langle \alpha
,\alpha \rangle $ so that $\langle \alpha ^{\vee },\alpha \rangle =2$. Also,
let $H_{\alpha }$ be defined by $\alpha \left( Z\right) =\langle H_{\alpha
},Z\rangle $, $Z\in \mathfrak{a}$, and write $H_{\alpha }^{\vee }=2H_{\alpha
}/\langle \alpha ,\alpha \rangle $ for the corresponding co-root. Finally,
let 
\[
\gamma _{\alpha }=\exp (i\pi H_{\alpha }^{\vee }), 
\]
where the exponential is taken in $\frak{g}_{\Bbb{C}}$, and put 
\[
F=\mbox{ group generated by }\{\gamma _{\alpha }:\alpha \in \Pi \}, 
\]
that is $F=\{\exp (i\pi H):H\in \mathcal{L}\}$, where $\mathcal{L}$ is the
lattice spanned by $H_{\alpha }^{\vee }$, $\alpha \in \Pi $.

It is known that $F$ is a subgroup of $M$ normalized by $M^{*}$ and that $%
M=FM_{0}$ (see Proposition 7.53 and Theorem 7.55 of \cite{knp}). Also, $%
\gamma _{\alpha }$ leaves invariant each root space $\frak{g}_{\beta }$ and
its restriction to $\frak{g}_{\beta }$ has the only eigenvalue $\exp ({i\pi
\langle \alpha ^{\vee },\beta \rangle })$. The next result shows that $F$
intersects each connected component of the centralizer $K_{H}$.

\begin{lema}
\label{lemacomponentes} For $H\in \mathfrak{a}$, we have that $%
K_{H}=F(K_{H})_{0}$. In particular, $K_{\Theta }=F(K_{\Theta })_{0}$.
\end{lema}

\begin{profe}
Take $w\in \mathcal{W}$ such that $Z=wH\in \mathrm{cl}\mathfrak{a}^{+}$.
Thus, since $K_{Z}=M(K_{Z})_{0}$ and $M=FM_{0}$, we have that $%
K_{Z}=F(K_{Z})_{0}$. Now 
\[
K_{H}=w^{-1}K_{Z}w=w^{-1}Fw(w^{-1}K_{Z}w)_{0}=F(K_{H})_{0}, 
\]
since $M^{*}$ normalizes $F$. The last assertion follows, since $K_{\Theta
}=K_{H_{\Theta }}$, where $H_{\Theta }\in \mathrm{cl}\frak{a}^{+}$ is such
that $\Theta (H_{\Theta })=\Theta $.
\end{profe}

\subsection{Stable and unstable bundles over the fixed points}
\label{secfibr}

Take $H\in \mathrm{cl}\frak{a}^{+}$. The one-parameter group $\exp (tH)$
acts on a flag manifold $\Bbb{F}_{\Theta }$, defining a flow, whose behavior
was described in Duistermat-Kolk-Varadarajan \cite{dkv}. This is the flow of
a gradient vector field, and the connected components of its fixed points
are given by the orbits $\mathrm{fix}_{\Theta }\left( H,w\right)
=K_{H}wb_{\Theta }$, where $w$ runs trough $\mathcal{W}$, $b_{\Theta }$ is
the origin of the flag manifold $\Bbb{F}_{\Theta }$ and $wb_{\Theta }=%
\overline{w}b_{\Theta }$, where $\overline{w}$ is any representative of $w$
in $M^{*}$. Since $K_{H}=K(H)M$ and the group $M$ fixes $wb_{\Theta }$, it
follows that 
\[
\mathrm{fix}_{\Theta }\left( H,w\right) =K(H)wb_{\Theta }. 
\]
It follows that $\mathrm{fix}_{\Theta }\left( H,w\right) =K\left( H\right)
/\left( K\left( H\right) \cap K_{wH_{\Theta }}\right) $, and hence $\mathrm{%
fix}_{\Theta }\left( H,w\right) $ is a flag manifold of the semisimple group 
$G(H)$.

The stable set of each $\mathrm{fix}_{\Theta }\left( H,w\right) $ is given
by 
\[
\mathrm{st}_{\Theta }(H,w)=N_{H}^{-}wb_{\Theta }, 
\]
and the stable bundle, denoted by $V_{\Theta }^{-}\left( H,w\right) $, is
the subbundle of the tangent bundle to $\mathrm{st}_{\Theta }(H,w)$
transversal to the fixed point set.

In order to write $V_{\Theta }^{-}\left( H,w\right) $ explicitly in terms of
root spaces we use the following notation: Given a vector subspace $\frak{l}%
\subset \frak{g}$ and $x\in \Bbb{F}_{\Theta }$ denote by $\frak{l}\cdot x$
the subspace of the tangent space $T_{x}\Bbb{F}_{\Theta }$ given by the
infinitesimal action of $\frak{l}$, namely 
\[
\frak{l}\cdot x=\{\widetilde{X}\left( x\right) \in T_{x}\Bbb{F}_{\Theta
}:X\in \frak{l}\}, 
\]
where $\widetilde{X}\left( x\right) =\frac{d}{dt}\left( \exp tX\right)
_{\left| t=0\right. }\left( x\right) $ is the vector field induced by $X\in 
\frak{g}$. With this notation the tangent space $T_{b_{\Theta }^{w}}\Bbb{F}%
_{\Theta }$ at $b_{\Theta }^{w}\approx wH_{\Theta }$ is 
\[
T_{b_{\Theta }^{w}}\Bbb{F}_{\Theta }=\frak{n}_{wH_{\Theta }}^{-}\cdot
b_{\Theta }^{w}. 
\]

Now, $V_{\Theta }^{-}\left( H,w\right) \rightarrow \mathrm{fix}_{\Theta
}\left( H,w\right) $ (which we write simpler as $V^{-}\rightarrow \mathrm{fix%
}_{\Theta }\left( H,w\right) $) is given by the following expressions:

\begin{enumerate}
\item  At $b_{\Theta }^{w}$ we put $V_{b_{\Theta }^{w}}^{-}=\left( \frak{n}%
_{wH_{\Theta }}^{-}\cap \frak{n}_{H}^{-}\right) \cdot b_{\Theta }^{w}$.

\item  At $x=gb_{\Theta }^{w}\in K_{H}\cdot b_{\Theta }^{w}$, $g\in K_{H}$
put 
\begin{equation}
V_{x}^{-}=\left( \mathrm{Ad}\left( g\right) \left( \frak{n}_{wH_{\Theta
}}^{-}\cap \frak{n}_{H}^{-}\right) \right) \cdot x.  \label{fordefve}
\end{equation}
This is the same as $dg_{b_{\Theta }^{w}}\left( V_{b_{\Theta }^{w}}\right) $
due to the well known formula $g_{*}\widetilde{X}=\widetilde{\left( \mathrm{%
Ad}\left( g\right) X\right) }$. Also, the right hand side of (\ref{fordefve}%
) depends only on $x$ because $\frak{n}_{wH_{\Theta }}^{-}\cap \frak{n}%
_{H}^{-}$ is invariant under the isotropy subgroup $K_{H}\cap K_{wH_{\Theta
}}$ of $\mathrm{fix}_{\Theta }\left( H,w\right) =K\left( H\right) /\left(
K\left( H\right) \cap K_{wH_{\Theta }}\right) $.
\end{enumerate}

For future reference we note that, by taking derivatives, the action of $%
K\left( H\right) $ on $\mathrm{fix}_{\Theta }\left( H,w\right) $ lifts to a
linear action on $V_{\Theta }^{-}\left( H,w\right) $. Also, in terms of root
spaces we have 
\[
\frak{n}_{wH_{\Theta }}^{-}\cap \frak{n}_{H}^{-}=\sum_{\beta \in \Pi
_{\Theta }^{-}\left( H,w\right) }\frak{g}_{\beta } 
\]
where 
\[
\Pi _{\Theta }^{-}\left( H,w\right) =\{\beta \in \Pi :\beta \left( H\right)
<0,\,\beta \left( wH_{\Theta }\right) <0\}. 
\]

In a similar way we can define the unstable bundles $V_{\Theta }^{+}\left(
H,w\right) \rightarrow \mathrm{fix}_{\Theta }\left( H,w\right) $ that are
tangent to the unstable sets $N_{H}^{+}wb_{\Theta }$ and transversal to the
fixed point set $\mathrm{fix}_{\Theta }\left( H,w\right) $. The construction
is the same unless that $\frak{n}_{H}^{-}$ is replaced by $\frak{n}_{H}^{+}$%
, and hence $\Pi _{\Theta }^{-}\left( H,w\right) $ is replaced by 
\[
\Pi _{\Theta }^{+}\left( H,w\right) =\{\beta \in \Pi :\beta \left( H\right)
>0,\,\beta \left( wH_{\Theta }\right) <0\}. 
\]

\vspace{12pt}%

\noindent \textbf{Remark:} The stable and unstable bundles $V_{\Theta }^{\pm
}\left( H,w\right) \rightarrow \mathrm{fix}_{\Theta }\left( H,w\right) $ can
be easily obtained by using the general device to construct a vector bundle
from a principal bundle $Q\rightarrow X$ and a representation of the
structural group $G$ on a vector space $V$. The resulting associated bundle $%
Q\times _{G}V$ is a vector bundle. For the stable and unstable bundles we
can take the principal bundle $K(H)\rightarrow \mathrm{fix}_{\Theta }\left(
H,w\right) $, defined by identification of $\mathrm{fix}_{\Theta }\left(
H,w\right) =K\left( H\right) /\left( K\left( H\right) \cap K_{wH_{\Theta
}}\right) $, whose structural group is $K\left( H\right) \cap K_{wH_{\Theta
}}$. Its representation on $\frak{l}^{\pm }=\frak{n}_{wH_{\Theta }}^{-}\cap 
\frak{n}_{H}^{\pm }$ yields $V_{\Theta }^{\pm }\left( H,w\right) $,
respectively.

\section{Vector bundles over homogeneous spaces}\label{secvechomog}

We state a general criterion of orientability of vector bundles acted by Lie
groups. Let $V\to X$ be a $n$-dimensional vector bundle and denote by $BV$
the bundle of frames $p:\Bbb{R}^{n}\rightarrow V$. It is well known that the
vector bundle $V$ is orientable if and only if $BV$ has exaclty two
connected components, and is connected otherwise.

Let $K$ be a connected Lie group acting transitively on the base space $X$
in such a way that the action lifts to a fiberwise linear action on $V$.
This linear action in turn lifts to an action on $BV$ by composition with
the frames.

Fix a base point $x_{0}\in X$ with isotropy subgroup $L\subset K$. Then each 
$g\in L$ gives rise to a linear operator of the fiber $\mathfrak{l}%
=V_{x_{0}} $. Denote by $\det (g|_{\mathfrak{l}})$, $g\in L$, the
determinant of this linear operator.

The following statement gives a simple criterion for the orientability of $V$%
.

\begin{proposicao}
\label{proporient}The vector bundle $V$ is orientable if and only if $\det
(g|_{\mathfrak{l}})>0$, for every $g\in L$.
\end{proposicao}

\begin{profe}
Suppose that $\det (g|_{\mathfrak{l}})>0$, $g\in L$, and take a basis $\beta
=\{e_{1},\ldots ,e_{k}\}$ of $V_{x_{0}}$. Let $g_{1},g_{2}\in G$ be such
that $g_{1}x_{0}=g_{2}x_{0}$. Then the bases $g_{i}\beta
=\{g_{i}e_{1},\ldots ,g_{i}e_{k}\}$, $i=1,2$, obtained by the linear action
on $V$, have the same orientation since $\deg \left( g_{1}^{-1}g_{2}|_{%
\mathfrak{l}}\right) >0$. These translations orient each fiber consistently
and hence $V$.

Conversely, denote by $BV$ the bundle of frames of $V$. If $V$ is orientable
then $BV$ splits into two connected components. Each one is a $\mathrm{Gl}%
^{+}\left( k,\Bbb{R}\right) $-subbundle, $k=\dim V$, and corresponds to an
orientation of $V$. The linear action of $G$ on $V$ lifts to an action on $%
BV $. Since $G$ is assumed to be connected, both connected components of $BV$
are $G$-invariant. Hence if $g\in L$ and $\beta $ is a basis of $V_{x_{0}}$
then $\beta $ and $g\beta $ have the same orientation, that is, $\det (g|_{%
\mathfrak{l}})>0$.
\end{profe}

\vspace{5pt} \noindent \textbf{Remark:} Clearly, $\det (g|_{\mathfrak{l}})$
does not change sign in a connected component of $L$. Hence to check the
condition of the above proposition it is enough to pick a point on each
connected component of $L$.

\subsection{Vector bundles over flag manifolds}\label{secorientaflag}

Now we are ready to get criteria for orientability of an stable vector
bundle $V_{\Theta }^{-}\left( H,w\right) \rightarrow \mathrm{fix}_{\Theta
}\left( H,w\right) $ and for the tangent bundle of a flag manifold $\Bbb{F}%
_{\Theta }$. These two cases have the following properties in common:

\begin{enumerate}
\item  The vector bundles are acted by a connected group whose action on the
base space is transitive. Hence Proposition \ref{proporient} applies.

\item  The connected components of the isotropy subgroup, at the base space,
is given by a subgroup $S$ of the lattice group $F$.

\item  The action of the isotropy subgroup on the fiber above the origin
reduces to the adjoint action on a space 
\[
\frak{l}=\sum_{\alpha \in \Gamma }\frak{g}_{\alpha } 
\]
spanned by root spaces, with roots belonging to a certain subset $\Gamma
\subset \Pi $.
\end{enumerate}

Now, a generator 
\[
\gamma _{\alpha }=\exp \left( i\pi H_{\alpha }^{\vee }\right) \qquad \alpha
\in \Pi 
\]
acts on a root space $\frak{g}_{\beta }$ by $\exp ({i\pi \langle \alpha
^{\vee },\beta \rangle })\cdot \mathrm{id}$. Hence the determinant of $%
\gamma _{\alpha }$ restricted to $\frak{l}=\sum_{\alpha \in \Gamma }\frak{g}%
_{\alpha }$ is given by 
\[
\det (\gamma _{\alpha }|_{\frak{l}})=\exp \left( i\pi \sum_{\beta \in \Gamma
}n_{\beta }\langle \alpha ^{\vee },\beta \rangle \right) . 
\]
So that $\det (\gamma _{\alpha }|_{\frak{l}})=\pm 1$ with the sign depending
whether the sum 
\[
\sum_{\beta \in \Gamma }n_{\beta }\langle \alpha ^{\vee },\beta \rangle 
\]
is even or odd. Here, as before $n_{\beta }$ is the multiplicity $\dim \frak{%
g}_{\beta }$ of the root $\beta $. From this we get the following criterion
for orientability in terms of roots: The vector bundle is orientable if and
only if for every root $\alpha $ the sum 
\[
\sum_{\beta \in \Gamma }n_{\beta }\langle \alpha ^{\vee },\beta \rangle
\equiv 0\qquad ({\rm mod}2) 
\]
where the sum is extended to $\beta \in \Gamma $.

\subsection{Flag manifolds}

In case of orientability of a flag manifold $\Bbb{F}_{\Theta }$ (its tangent
bundle) the subspace to be considered is 
\[
\frak{l}=\frak{n}_{{\Theta }}^{-}=\sum_{\beta \in \Pi ^{-}\setminus \langle
\Theta \rangle }\frak{g}_{\beta }, 
\]
that identifies with the tangent space to $\Bbb{F}_{\Theta }$ at the origin.
On the other hand the isotropy subgroup $K_{\Theta }=F(K_{\Theta })_{0}$
(see Lemma \ref{lemacomponentes}), which means that $F$ covers the connected
components of $K_{\Theta }$. Hence we get the following criterion.

\begin{teorema}
\label{teoorientflags}The flag manifold $\Bbb{F}_{\Theta }$ is orientable if
and only if 
\begin{equation}
\sum_{\beta }n_{\beta }\langle \alpha ^{\vee },\beta \rangle \equiv 0\quad (%
\mathrm{mod}2)  \label{fororientflags}
\end{equation}
where the sum is extended to $\beta \in \Pi ^{-}\setminus \langle \Theta
\rangle $ (or equivalently to $\beta \in \Pi ^{+}\setminus \langle \Theta
\rangle $). This condition must be satisfied for any simple root $\alpha $.
\end{teorema}

\begin{profe}
In fact, $\Pi ^{-}\setminus \langle \Theta \rangle $ is the set of roots
whose root spaces span the tangent space at the origin. Hence the
determinant condition holds if (\ref{fororientflags}) is satisfied for every
root $\alpha \in \Pi $. However it is enough to take $\alpha $ in the simple
system $\Sigma $. This is because the set of co-roots $\Pi ^{\vee }=\{\alpha
^{\vee }:\alpha \in \Pi \}$ is also a root system having $\Sigma ^{\vee
}=\{\alpha ^{\vee }:\alpha \in \Sigma \}$ as a simple system of roots. By
taking linear combinations of $\Sigma ^{\vee }$ with integer coefficients it
follows that condition (\ref{fororientflags}) holds for any root $\alpha \in
\Pi $ if and only if it is satisfied for the simple roots.
\end{profe}

Now we derive some consequences of the criteria stated above. First we prove
that any maximal flag manifold is orientable, a result already obtained by
Kocherlakota \cite{koch} as a consequence that the top $\Bbb{Z}$-homology
groups are nontrivial.

\begin{teorema}
Any maximal flag manifold $\Bbb{F}$ is orientable.
\end{teorema}

\begin{profe}
We write, for a simple root $\alpha $, $\Pi _{\alpha }=\{\alpha ,2\alpha
\}\cap \Pi ^{+}$, $\Pi _{0}^{\alpha }=\{\beta \in \Pi ^{+}:\langle \alpha
^{\vee },\beta \rangle =0\}$ and $\Pi _{1}^{\alpha }=\{\beta \in \Pi
^{+}:\langle \alpha ^{\vee },\beta \rangle \neq 0$, $\beta \notin \Pi
_{\alpha }\}$. Let $r_{\alpha }$ be the reflection with respect to $\alpha $%
. It is known that $r_{\alpha }\left( \Pi ^{+}\setminus \Pi _{\alpha
}\right) =\Pi ^{+}\setminus \Pi _{\alpha }$. Moreover, for a root $\beta $
we have 
\[
\langle \alpha ^{\vee },r_{\alpha }\left( \beta \right) \rangle =\langle
\alpha ^{\vee },\beta -\langle \alpha ^{\vee },\beta \rangle \alpha \rangle
=\langle \alpha ^{\vee },\beta \rangle -\langle \alpha ^{\vee },\alpha
\rangle \langle \alpha ^{\vee },\beta \rangle =-\langle \alpha ^{\vee
},\beta \rangle . 
\]
Hence the subsets $\Pi _{0}^{\alpha }$ and $\Pi _{1}^{\alpha }$ are $%
r_{\alpha }$-invariant and $\langle \alpha ^{\vee },\beta +r_{\alpha }\left(
\beta \right) \rangle =0$.

Now fix $\alpha \in \Sigma $ and split the sum $\sum_{\beta \in \Pi
^{+}}n_{\beta }\langle \alpha ^{\vee },\beta \rangle $ into $\Pi _{\alpha }$%
, $\Pi _{0}^{\alpha }$ and $\Pi _{1}^{\alpha }$. For $\Pi _{\alpha }$ this
sum is $2n_{\alpha }+4n_{2\alpha }$, with $n_{2\alpha }=0$ if $2\alpha $ is
not a root. For $\Pi _{0}^{\alpha }$ the sum is zero. In $\Pi _{1}^{\alpha }$
the roots are given in pairs $\beta \neq r_{\alpha }\left( \beta \right) $
with $\langle \alpha ^{\vee },\beta +r_{\alpha }\left( \beta \right) \rangle
=0$, since $\Pi _{1}^{\alpha }$ is $r_{\alpha }$-invariant and $\beta
=r_{\alpha }\left( \beta \right) $ if and only if $\langle \alpha ^{\vee
},\beta \rangle =0$. Since $n_{r_{\alpha }\left( \beta \right) }=n_{\beta }$%
, it follows that $\sum_{\beta \in \Pi _{1}^{\alpha }}n_{\beta }\langle
\alpha ^{\vee },\beta \rangle =0$. Hence the total sum is even for every $%
\alpha \in \Sigma $, proving the orientability of $\Bbb{F}$.
\end{profe}

In particular this orientability result applies to the maximal flag manifold
of the semi-simple Lie algebra $\frak{g}(\Theta )$. Here the set of roots is 
$\langle \Theta \rangle $ having $\Theta $ as a simple system of roots.
Therefore the equivalent conditions of Theorem \ref{teoorientflags} combined
with the orientability of the maximal flag manifold of $\frak{g}(\Theta )$
immplies the

\begin{corolario}
\label{corsomaflagTheta}If $\alpha \in \Theta $ then 
\[
\sum_{\beta }n_{\beta }\langle \alpha ^{\vee },\beta \rangle \equiv 0\quad (%
\mathrm{mod}2), 
\]
where the sum is extended to $\beta \in \langle \Theta \rangle ^{-}$ (or
equivalently to $\beta \in \langle \Theta \rangle ^{+}$).
\end{corolario}

This allows to simplify the criterion for a partial flag manifold $\Bbb{F}%
_{\Theta }$.

\begin{proposicao}
\label{propsomaflagparcial}$\Bbb{F}_{\Theta }$ is orientable if and only if,
for every root $\alpha \in \Sigma \setminus \Theta $, it holds 
\begin{equation}
\sum_{\beta }n_{\beta }\langle \alpha ^{\vee },\beta \rangle \equiv 0\quad (%
\mathrm{mod}2),  \label{forsomaflagparcial}
\end{equation}
where the sum is extended to $\beta \in \langle \Theta \rangle ^{-}$ (or
equivalently to $\beta \in \langle \Theta \rangle ^{+}$).
\end{proposicao}

\begin{profe}
Applying Corollary \ref{corsomaflagTheta} with $\Theta =\Sigma $, we have
that $\sum_{\beta \in \Pi ^{-}}n_{\beta }\langle \alpha ^{\vee },\beta
\rangle $ is even. Hence, by Theorem \ref{teoorientflags}, $\Bbb{F}_{\Theta
} $ is orientable if and only if, for every root $\alpha \in \Sigma $, the
sum $\sum_{\beta \in \langle \Theta \rangle ^{-}}n_{\beta }\langle \alpha
^{\vee },\beta \rangle $ is even. By Corollary \ref{corsomaflagTheta}, it is
enough to check this for every root $\alpha \in \Sigma \setminus \Theta $.
\end{profe}

Finally we observe that if $G$ is a complex group then the real
multiplicities are $n_{\beta }=2$ so that any flag $\Bbb{F}_{\Theta }$ is
orientable. This is well known since the $\Bbb{F}_{\Theta }$ are complex
manifolds.

\subsection{Stable and unstable bundles in flag manifolds}

For the stable bundles $V_{\Theta }^{-}\left( H,w\right) $ we take 
\[
\frak{l}=\frak{n}_{wH_{\Theta }}^{-}\cap \frak{n}_{H}^{-}=\sum_{\beta \in
\Pi _{\Theta }^{-}\left( H,w\right) }\frak{g}_{\beta }. 
\]
where 
\[
\Pi _{\Theta }^{-}\left( H,w\right) =\{\beta \in \Pi :\beta \left( H\right)
<0,\,\beta \left( wH_{\Theta }\right) <0\}. 
\]
Also the acting Lie group is $K(H)$ whose isotropy subgroup at $wH_{\Theta }$
of the base space $\mathrm{fix}_{\Theta }\left( H,w\right) $ is $L=K(H)\cap
Z_{wH_{\Theta }}$ where $Z_{wH_{\Theta }}$ is the centralizer of $wH_{\Theta
}$. Applying the determinant criterion we get the following condition for
orientability.

\begin{teorema}
\label{teoorientstable}The vector bundle $V_{\Theta }^{-}\left( H,w\right) $
is orientable if and only if 
\[
\sum_{\beta }n_{\beta }\langle \alpha ^{\vee },\beta \rangle \equiv 0\quad (%
\mathrm{mod}2), 
\]
where the sum is extended to $\beta \in \Pi _{\Theta }^{-}\left( H,w\right) $%
. Here the condition must be verified for every $\alpha \in \Theta (H)$.
\end{teorema}

\begin{profe}
It remains to discuss the last statement about the scope of the condition.
It is a consequence of Lemma \ref{lemacomponentes}. In fact, $K(H)$ is the
compact component of the semisimple Lie group $G(H)$. Hence 
\[
L=K(H)\cap Z_{wH_{\Theta }}=F(H)(K(H)\cap Z_{wH_{\Theta }})_{0}, 
\]
where $F(H)$ is the $F$ group of $G(H)$, that is, the group generated by 
\[
\{\gamma _{\alpha }=\exp \left( i\pi H_{\alpha }^{\vee }\right) :\,\alpha
\in \langle \Theta (H)\rangle \}, 
\]
because the restriction of $\langle \Theta (H)\rangle $ to $\frak{a}(H)$ is
the root system of $G(H)$. Finally, it is enough to check the condition for
the simple roots in $\Theta \left( H\right) $.
\end{profe}

\vspace{12pt}%

\noindent \textbf{Remark:} The same result holds for the unstable vector
bundles $V_{\Theta }^{+}\left( H,w\right) $ with $\Pi _{\Theta }^{+}\left(
H,w\right) $ instead of $\Pi _{\Theta }^{-}\left( H,w\right) $.

We have the following result in the special case when $\Theta =\emptyset $
and $w$ is the principal involution $w^{-}$.

\begin{corolario}
For every $H\in \mathrm{cl}\frak{a}^{+}$, the vector bundles $V^{-}\left(
H,1\right) $ and $V^{+}\left( H,w^{-}\right) $ are orientable.
\end{corolario}

\begin{profe}
Applying Corollary \ref{corsomaflagTheta} with $\Theta =\Sigma $ and $\Theta
=\Theta (H)$, it follows that both 
\[
\sum_{\beta \in \Pi ^{+}}n_{\beta }\langle \alpha ^{\vee },\beta \rangle
\qquad \mbox{and}\qquad \sum_{\beta \in \langle \Theta (H)\rangle
^{+}}n_{\beta }\langle \alpha ^{\vee },\beta \rangle 
\]
are even for $\alpha \in \Theta (H)$. Hence, for every $\alpha \in \Theta
(H) $, it holds that $\sum_{\beta }n_{\beta }\langle \alpha ^{\vee },\beta
\rangle $ is even, where the sum is extended to $\beta \in \Pi ^{+}\setminus
\langle \Theta (H)\rangle $. If $\Theta =\emptyset $, then $H_{\Theta }$ is
regular and $\beta \left( w^{-}H_{\Theta }\right) <0$ if and only if $\beta
\in \Pi ^{+}$. Thus $\Pi ^{+}\left( H,w^{-}\right) =\Pi ^{+}\setminus
\langle \Theta (H)\rangle $ and the result follows from Theorem \ref
{teoorientstable}.

The proof for $V^{+}\left( H,w^{-}\right) $ is analogous.
\end{profe}

\vspace{12pt}%

\noindent \textbf{Remark:} The above result is not true in a partial flag
manifold. An example is given in $G=\mathrm{Sl}\left( 3,\Bbb{R}\right) $
with $H=\mathrm{diag}\{2,-1,-1\}$. Then it can be seen that the repeller
component of $H$ is a projective line and its unstable bundle a M\"{o}bius
strip.

\section{Split real forms}

\label{secsplit}

When $\frak{g\,}$is a split real form every root $\beta $ has multiplicity $%
n_{\beta }=1$. Hence, the criterion of Corollary \ref{propsomaflagparcial}
reduces to 
\begin{equation}
S\left( \alpha ,\Theta \right) =\sum_{\beta \in \langle \Theta \rangle
^{+}}\langle \alpha ^{\vee },\beta \rangle \equiv 0\qquad \left( \mathrm{mod}%
2\right) ,  \label{forsomasemmult}
\end{equation}
that can be checked by looking at the Dynkin diagrams. In the sequel we use
a standard way of labelling the roots in the diagrams as in the picture
below.

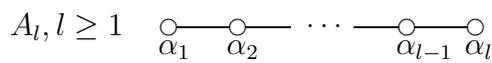
\begin{picture}(400,330)(-20,0)

\put(-20,260){
\begin{picture}(140,60)(-40,0)

\put(-40,27){$A_l,l \geq 1$}
\put(20,30){\circle{6}}
\put(23,30){\line(1,0){20}}
\put(46,30){\circle{6}}
\put(49,30){\line(1,0){17}}
\put(72,30){\ldots}
\put(90,30){\line(1,0){17}}
\put(110,30){\circle{6}}
\put(113,30){\line(1,0){20}}
\put(136,30){\circle{6}}
\put(16,20){${\alpha}_1$}
\put(42,20){${\alpha}_2$}
\put(106,20){${\alpha}_{l-1}$}
\put(132,20){${\alpha}_l$}
\end{picture}
}

\put(-20,200)
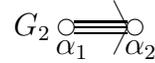
{
\begin{picture}(100,60)(-40,0)

\put(-40,27){$B_l,l \geq 2$}
\put(20,30){\circle{6}}
\put(23,30){\line(1,0){20}}
\put(46,30){\circle{6}}
\put(49,30){\line(1,0){17}}
\put(72,30){\ldots}
\put(90,30){\line(1,0){17}}
\put(110,30){\circle{6}}
\put(113,31.2){\line(1,0){20}}
\put(113,28.8){\line(1,0){20}}
\put(136,30){\circle{6}}
\put(133,30){\line(-1,2){5}}
\put(133,30){\line(-1,-2){5}}
\put(16,20){${\alpha}_1$}
\put(42,20){${\alpha}_2$}
\put(106,20){${\alpha}_{l-1}$}
\put(132,20){${\alpha}_l$}

\end{picture}
}
\put(-20,140)
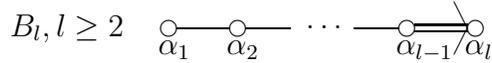
{
\begin{picture}(100,60)(-40,0)

\put(-40,27){$C_l,l \geq 3$}
\put(20,30){\circle{6}}
\put(23,30){\line(1,0){20}}
\put(46,30){\circle{6}}
\put(49,30){\line(1,0){17}}
\put(72,30){\ldots}
\put(90,30){\line(1,0){17}}
\put(110,30){\circle{6}}
\put(113,30){\line(1,2){5}}
\put(113,30){\line(1,-2){5}}
\put(113,31.2){\line(1,0){20}}
\put(113,28.8){\line(1,0){20}}
\put(136,30){\circle{6}}
\put(16,20){${\alpha}_1$}
\put(42,20){${\alpha}_2$}
\put(106,20){${\alpha}_{l-1}$}
\put(132,20){${\alpha}_l$}
\end{picture}
}
\put(-20,70)
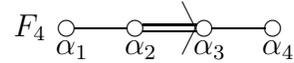
{
\begin{picture}(100,70)(-40,0)

\put(-40,32){$D_l,l \geq 4$}
\put(20,35){\circle{6}}
\put(16,25){$\alpha_1$}
\put(23,35){\line(1,0){20}}
\put(46,35){\circle{6}}
\put(42,25){$\alpha_2$}
\put(49,35){\line(1,0){17}}
\put(72,35){\ldots}
\put(90,35){\line(1,0){17}}
\put(110,35){\circle{6}}
\put(105,25){$\alpha_{l-2}$}
\put(111.9,36.6){\line(6,5){20}}
\put(111.9,33.4){\line(6,-5){20}}
\put(134,54.6){\circle{6}}
\put(139,51){$\alpha_{l-1}$}
\put(134,15.6){\circle{6}}
\put(139,12){$\alpha_l$}

\end{picture}
}
\put(220,260)
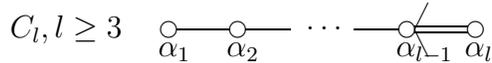
{
\begin{picture}(100,60)(0,0)

\put(0,27){$G_2$}
\put(20,30){\circle{6}}
\put(23,30){\line(1,0){20}}
\put(23,32.2){\line(1,0){20}}
\put(23,27.8){\line(1,0){20}}
\put(46,30){\circle{6}}
\put(43,30){\line(-1,2){5}}
\put(43,30){\line(-1,-2){5}}
\put(16,20){${\alpha}_1$}
\put(42,20){${\alpha}_2$}
\end{picture}
}
\put(220,200)
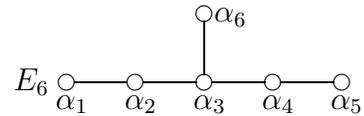
{
\begin{picture}(100,60)(0,0)

\put(0,27){$F_4$}
\put(20,30){\circle{6}}
\put(16,20){${\alpha}_1$}
\put(23,30){\line(1,0){20}}
\put(46,30){\circle{6}}
\put(42,20){${\alpha}_2$}
\put(49,31.2){\line(1,0){20}}
\put(49,28.8){\line(1,0){20}}
\put(72,30){\circle{6}}
\put(68,20){${\alpha}_3$}
\put(69,30){\line(-1,2){5}}
\put(69,30){\line(-1,-2){5}}
\put(75,30){\line(1,0){20}}
\put(98,30){\circle{6}}
\put(94,20){${\alpha}_4$}
\end{picture}
}
\put(220,128)
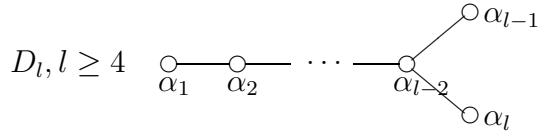
{
\begin{picture}(130,70)(0,0)

\put(0,27){$E_6$}
\put(20,30){\circle{6}}
\put(23,30){\line(1,0){20}}
\put(46,30){\circle{6}}
\put(49,30){\line(1,0){20}}
\put(72,30){\circle{6}}
\put(75,30){\line(1,0){20}}
\put(98,30){\circle{6}}
\put(101,30){\line(1,0){20}}
\put(124,30){\circle{6}}
\put(72,33){\line(0,1){20}}
\put(72,56){\circle{6}}
\put(16,20){${\alpha}_1$}
\put(42,20){${\alpha}_2$}
\put(68,20){${\alpha}_3$}
\put(94,20){${\alpha}_4$}
\put(120,20){${\alpha}_5$}
\put(76,53){${\alpha}_6$}

\end{picture}
}
\put(220,66)
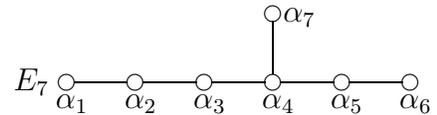
{
\begin{picture}(200,70)(0,0)

\put(0,27){$E_7$}
\put(20,30){\circle{6}}
\put(23,30){\line(1,0){20}}
\put(46,30){\circle{6}}
\put(49,30){\line(1,0){20}}
\put(72,30){\circle{6}}
\put(75,30){\line(1,0){20}}
\put(98,30){\circle{6}}
\put(101,30){\line(1,0){20}}
\put(124,30){\circle{6}}
\put(127,30){\line(1,0){20}}
\put(150,30){\circle{6}}
\put(98,33){\line(0,1){20}}
\put(98,56){\circle{6}}
\put(16,20){${\alpha}_1$}
\put(42,20){${\alpha}_2$}
\put(68,20){${\alpha}_3$}
\put(94,20){${\alpha}_4$}
\put(120,20){${\alpha}_5$}
\put(146,20){${\alpha}_6$}
\put(102,53){${\alpha}_7$}

\end{picture}
}
\put(220,0)
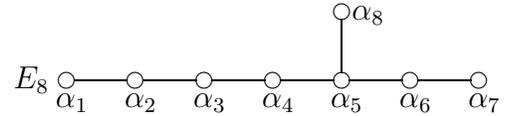
{
\begin{picture}(200,70)(0,0)

\put(0,27){$E_8$}
\put(20,30){\circle{6}}
\put(23,30){\line(1,0){20}}
\put(46,30){\circle{6}}
\put(49,30){\line(1,0){20}}
\put(72,30){\circle{6}}
\put(75,30){\line(1,0){20}}
\put(98,30){\circle{6}}
\put(101,30){\line(1,0){20}}
\put(124,30){\circle{6}}
\put(127,30){\line(1,0){20}}
\put(150,30){\circle{6}}
\put(153,30){\line(1,0){20}}
\put(176,30){\circle{6}}
\put(124,33){\line(0,1){20}}
\put(124,56){\circle{6}}
\put(16,20){${\alpha}_1$}
\put(42,20){${\alpha}_2$}
\put(68,20){${\alpha}_3$}
\put(94,20){${\alpha}_4$}
\put(120,20){${\alpha}_5$}
\put(146,20){${\alpha}_6$}
\put(172,20){${\alpha}_7$}
\put(128,53){${\alpha}_8$}

\end{picture}
}

\end{picture}

For the diagram $G_{2}$ there are three flag manifolds: the maximal $\Bbb{F}$%
, which is orientable, and the minimal ones $\Bbb{F}_{\{\alpha _{1}\}}$ and $%
\Bbb{F}_{\{\alpha _{2}\}}$, where $\alpha _{1}$ and $\alpha _{2}$ are the
simple roots with $\alpha _{1}$ the longer one. These minimal flag manifolds
are not orientable since in both cases (\ref{forsomasemmult}) reduces to the
Killing numbers $\langle \alpha _{1}^{\vee },\alpha _{2}\rangle =-1$ and $%
\langle \alpha _{2}^{\vee },\alpha _{1}\rangle =-3$. From now on we consider
only simple and double laced diagrams.

Our strategy consists in counting the contribution of each connected
component $\Delta $ of $\Theta $ to the sum $S\left( \alpha ,\Theta \right) $
in (\ref{forsomasemmult}). Thus we keep fixed $\alpha $ and a connected
subset $\Delta \subset \Sigma $. If $\alpha $ is not linked to $\Delta $
then $S\left( \alpha ,\Delta \right) =0$ and we can discard this case.
Otherwise, $\alpha $ is linked to exactly one root of $\Delta $, because a
Dynkin diagram has no cycles. We denote by $\delta $ the only root in $%
\Delta $ linked to $\alpha $.

A glance at the Dynkin diagrams show the possible subdiagrams $\Delta $
properly contained in $\Sigma $. We exhibit them in table \ref{tabeltodos}.
For these subdiagrams we can write down explicitly the roots of $\langle
\Delta \rangle ^{+}$ and then compute $S\left( \alpha ,\Delta \right) $,
when $\alpha $ is linked to $\Delta $. In fact, if $\beta \subset \langle
\Delta \rangle ^{+}$ then $\beta =c\delta +\gamma $ where $\delta $ is the
only root in $\Delta $ which is linked to $\alpha $ and $\langle \gamma
,\alpha ^{\vee }\rangle =0$, so that $\langle \beta ,\alpha ^{\vee }\rangle
=c\langle \delta ,\alpha ^{\vee }\rangle $. Hence it is enough to look at
those roots $\beta \in \Delta $ whose coefficient $c$ in the direction of $%
\delta $ is nonzero. In the sequel we write down the values of $S\left(
\alpha ,\Delta \right) $ and explain how they were obtained.

\begin{table}[tbp]
\centering
\par
\begin{tabular}{|l|c|}
\hline
\,\,\,\,\,\,\,\,\,\, $\Delta $ & $\Sigma $ \\ \hline
$A_k$ ($k \geq 1 $) & any diagram \\ \hline
$B_k $ ($k \geq 2 $) & $B_l $ ($l>k$), $C_l $ ($k = 2$) and $F_4 $ ($2 \leq
k \leq 3$) \\ \hline
$C_k $ ($k \geq 3 $) & $C_l $ ($l>k$) and $F_4 $ ($k = 3$) \\ \hline
$D_k $ ($k \geq 4 $) & $D_l $ ($l>k$), $E_6 $ ($4 \leq k \leq 5$), $E_7 $ ($%
4 \leq k \leq 6$) and $E_8 $ ($4 \leq k \leq 7$) \\ \hline
$E_6$ & $E_7$ and $E_8$ \\ \hline
$E_7$ & $E_8$ \\ \hline
\end{tabular}
\caption{Connected subdiagrams}
\label{tabeltodos}
\end{table}

In the diagram $A_{k}$ with roots $\alpha _{1},\ldots ,\alpha _{k}$ the
positive roots are $\alpha _{i}+\cdots +\alpha _{j}$, $i\leq j$. Hence if $%
\Delta =A_{k}$ then the possibilities for $\delta $ are the extreme roots $%
\alpha _{1}$ and $\alpha _{k}$. In case $\delta =\alpha _{1}$ the sum $%
S\left( \alpha ,\Delta \right) $ extends over the $k$ positive roots $\alpha
_{1}+\cdots +\alpha _{j}$, $j=1,\ldots ,k$, that have nonzero coefficient in
the direction of $\alpha _{1}$. (It is analogous for $\delta =\alpha _{k}$.)

\begin{table}[tbp]
\centering
\par
\begin{tabular}{|l|c|}
\hline
\multicolumn{2}{|c|}{$\Delta = A_k $} \\ \hline
links & $S(\alpha , \Delta )$ \\ \hline
\begin{picture}(40,15)(0,0) \put(0,0){$\alpha $} \put(10,3){\line(1,0){20}}
\put(33,0){$\delta $} \end{picture} & $- k$ \\ \hline
\begin{picture}(40,15)(0,0) \put(0,3){$\alpha $} \put(30,3){\line(-1,2){5}}
\put(30,6){\line(-1,-2){5}} \put(10,3){\line(1,0){20}}
\put(10,6){\line(1,0){20}} \put(33,3){$\delta $} \end{picture} & $- k$ \\ 
\hline
\begin{picture}(40,15)(0,0) \put(0,3){$\alpha $} \put(10,3){\line(1,2){5}}
\put(10,6){\line(1,-2){5}} \put(10,3){\line(1,0){20}}
\put(10,6){\line(1,0){20}} \put(33,3){$\delta $} \end{picture} & $- 2k$ \\ 
\hline
\end{tabular}
\caption{$A_l$ subdiagrams}
\label{tabaele}
\end{table}

In the standard realization of $B_{k}$ the positive roots are $\lambda
_{i}\pm \lambda _{j}$, $i\neq j$, and $\lambda _{i}$, where $\{\lambda
_{1},\ldots ,\lambda _{k}\}$ is an orthonormal basis of the $k$-dimensional
space. The possibilities for $\delta $ are extreme roots $\lambda
_{1}-\lambda _{2}$ (to the left) and $\lambda _{k}$ (to the right). If $%
\delta =\lambda _{1}-\lambda _{2}$ then $\alpha $ and $\delta $ are linked
by one edge, that is, $\langle \delta ,\alpha ^{\vee }\rangle =-1$. Also,
the positive roots in $B_{k}$ having nonzero coefficient $c$ in the
direction of $\lambda _{1}-\lambda _{2}$ are the $2k-2$ roots $\lambda
_{1}\pm \lambda _{j}$, $j>1$ together with $\lambda _{1}$. For all of them $%
c=1$, hence the contribution of $\Delta $ to $S\left( \alpha ,\Delta \right) 
$ is $-\left( 2k-1\right) $. Analogous computations with $\delta =\lambda
_{k}$ yields the table

\begin{table}[tbp]
\centering
\par
\begin{tabular}{c}
$\Delta = B_k $ \\ 
\begin{tabular}{|l|c|}
\hline
$\Sigma$ & $S(\alpha , \Delta )$ \\ \hline
$B_l $ ($2 \leq k < l$) & $- (2k - 1)$ \\ \hline
$C_l $ ($k = 2$) & $-4$ \\ \hline
$F_4 $ ($k = 2$) & $-3$ or $-4$ \\ \hline
$F_4 $ ($k = 3$) & $-9$ \\ \hline
\end{tabular}
\end{tabular}
\caption{$B_l$ subdiagrams}
\label{tabbele}
\end{table}

For $C_{k}$ the positive roots are $\lambda _{i}\pm \lambda _{j}$, $i\neq j$%
, and $2\lambda _{i}$. If $\delta =\lambda _{1}-\lambda _{2}$ then $\langle
\delta ,\alpha ^{\vee }\rangle =-1$, and we must count the $2k-2$ roots $%
\lambda _{1}\pm \lambda _{j}$, $j>1$, having coefficient $c=1$ and $2\lambda
_{1}$ with $c=2$. Then the contribution to $S\left( \alpha ,\Delta \right) $
is $-2k$. This together with a similar computation for the other $\delta $
gives table

\begin{table}[tbp]
\centering
\par
\begin{tabular}{c}
$\Delta = C_k $ \\ 
\begin{tabular}{|l|c|}
\hline
$\Sigma$ & $S(\alpha , \Delta )$ \\ \hline
$C_l $ ($3 \leq k < l$) & $- 2k $ \\ \hline
$F_4 $ ($k = 3$) & $-6$ \\ \hline
\end{tabular}
\end{tabular}
\caption{$C_l$ subdiagrams}
\label{tabcele}
\end{table}

For $D_{k}$ the positive roots are $\lambda _{i}\pm \lambda _{j}$, $i\neq j$%
. If $\delta =\lambda _{1}-\lambda _{2}$ then $\langle \delta ,\alpha ^{\vee
}\rangle =-1$, and we must count the $2k-2$ roots $\lambda _{1}\pm \lambda
_{j}$, $j>1$, all of them having coefficient $c=1$. Then the contribution to 
$S\left( \alpha ,\Delta \right) $ is $-2k-2$. We leave to the reader the
computation of the other entries of table

\begin{table}[tbp]
\centering
\par
\begin{tabular}{c}
$\Delta = D_k $ \\ 
\begin{tabular}{|l|c|}
\hline
$\Sigma$ & $S(\alpha , \Delta )$ \\ \hline
$D_l $ ($4 \leq k < l$) & $- 2(k - 1)$ \\ \hline
$E_l $ ($k = 4$) & $-6$ \\ \hline
$E_l $ ($k = 5$) & $- 8$, $\delta = \alpha _1$ \\ \hline
$E_l $ ($k = 5$) & $- 10$, $\delta = \alpha _5$ \\ \hline
$E_l $ ($k = 6$) & $- 6$, $\delta = \alpha _1$ \\ \hline
$E_l $ ($k = 6$) & $- 15$, $\delta = \alpha _6$ \\ \hline
$E_8 $ ($k = 7$) & $-21$ \\ \hline
\end{tabular}
\end{tabular}
\caption{$D_l$ subdiagrams}
\label{tabdele}
\end{table}

The results for the exceptional cases are included in table \ref{tabele}. To
do the computations we used the realization of Freudenthal of the split real
form of $E_{8}$ in the vector space $\frak{sl}\left( 9,\Bbb{R}\right) \oplus
\bigwedge^{3}\Bbb{R}^{9}\oplus \left( \bigwedge^{3}\Bbb{R}^{9}\right) ^{*}$.
The roots of $E_{8}$ are the weights of the representation of the Cartan
subalgebra $\frak{h}\subset \frak{sl}\left( 9,\Bbb{R}\right) $ of the
diagonal matrices (see Fulton-Harris \cite{fuha} and \cite{smalg}). The
roots are $\lambda _{i}-\lambda _{j}$, $i\neq j$ (with root spaces in $\frak{%
sl}\left( 9,\Bbb{R}\right) $) and $\pm \left( \lambda _{i}+\lambda
_{j}+\lambda _{k}\right) $, $i<j<k$ (with root spaces in $\bigwedge^{3}\Bbb{R%
}^{9}\oplus \left( \bigwedge^{3}\Bbb{R}^{9}\right) ^{*}$). From the
realization of $E_{8}$ one easily obtains $E_{6}$ and $E_{7}$, and the
computations can be performed.

\begin{table}[tbp]
\centering
\par
\begin{tabular}{c}
$\Delta = E_k $ \\ 
\begin{tabular}{|l|c|}
\hline
$\Sigma$ & $S(\alpha , \Delta )$ \\ \hline
$E_l $ ($k = 6$) & $- 16 $ \\ \hline
$E_8 $ ($k = 7$) & $-27 $ \\ \hline
\end{tabular}
\end{tabular}
\caption{$E_l$ subdiagrams}
\label{tabele}
\end{table}

\subsection{Classical Lie algebras}

The split real forms of the classical Lie algebras are $A_{l}=\frak{sl}%
\left( l+1,\Bbb{R}\right) $, $B_{l}=\frak{so}\left( l,l+1\right) $, $C_{l}=%
\frak{sp}\left( l,\Bbb{R}\right) $ and $D_{l}=\frak{so}\left( l,l\right) $.
Their associated flag manifolds are concretely realized as manifolds of
flags $\left( V_{1}\subset \cdots \subset V_{k}\right) $ of vector subspaces 
$V_{i}\subset \Bbb{R}^{n}$. For $A_{l}$ one take arbitrary subspaces of $%
\Bbb{R}^{n}$, $n=l+1$. Given integers $1\leq d_{1}<\cdots <d_{k}\leq l$ we
denote by $\Bbb{F}\left( d_{1},\ldots ,d_{k}\right) $ the manifold of flags $%
\left( V_{1}\subset \cdots \subset V_{k}\right) $ with $\dim V_{i}=d_{i}$.

For the other classical Lie algebras we take similar manifolds of flags, but
now the subspaces $V_{i}$ are isotropic w.r.t. a quadratic form for $B_{l}$
and $D_{l}$, and w.r.t. a symplectic form in $C_{l}$. Again the flag
manifolds are given by integers $1\leq d_{1}<\cdots <d_{k}\leq l$ and we
write $\Bbb{F}^{I}\left( d_{1},\ldots ,d_{k}\right) $ for the manifold of
flags of isotropic subspaces with $\dim V_{i}=d_{i}$. Here $V_{i}\subset 
\Bbb{R}^{n}$ with $n=2l+1$ in $B_{l}$ and $n=2l$ in the $C_{l}$ and $D_{l}$
cases.

The way we order the simple roots $\Sigma $ in the Dynkin diagrams allows a
direct transition between the dimensions $d_{1},\ldots ,d_{k}$ and the roots 
$\Theta \subset \Sigma $ when $\Bbb{F}\left( d_{1},\ldots ,d_{k}\right) $ or 
$\Bbb{F}^{I}\left( d_{1},\ldots ,d_{k}\right) $ is $\Bbb{F}_{\Theta }$. In
fact, except for some flags of $D_{l}$ the dimensions $d_{1},\ldots ,d_{k}$
coincide with the indices of the roots $\alpha _{j}\notin \Theta $. (For
example, the Grassmannian $\Bbb{F}\left( r\right) $ is the flag manifold $%
\Bbb{F}_{\Theta }$ with $\Theta =\Sigma \setminus \{\alpha _{r}\}$.) We
detail this correspondence below.

The orientability criteria for the split real groups uses several times the
following

\vspace{12pt}%

\noindent \textbf{Condition:} We say that the numbers $0=d_{0},d_{1},\ldots
,d_{k}$ satisfy the ${\rm mod}2$ condition if the differences $%
d_{i+1}-d_{i} $, $i=0,\ldots ,k$, are congruent ${\rm mod}2$, that is, they
are simultaneously even or simultaneously odd.

\subsubsection{$A_{l}=\frak{sl}\left( l+1,\Bbb{R}\right) $}

The flag manifolds are $\Bbb{F}\left( d_{1},\ldots ,d_{k}\right) =\Bbb{F}%
_{\Theta }$ such that $j\in \{d_{1},\ldots ,d_{k}\}$ if and only if $j$ is
the index of a simple root $\alpha _{j}\notin \Theta $. If we write $\Bbb{F}%
\left( d_{1},\ldots ,d_{k}\right) =\mathrm{SO}\left( n\right) /K_{\Theta }$
then $K_{\Theta }=\mathrm{SO}\left( d_{1}\right) \times \cdots \times 
\mathrm{SO}\left( n-d_{k}\right) $ is a group of block diagonal matrices,
having blocks of sizes $d_{i+1}-d_{i}$.

\begin{proposicao}
A flag manifold $\Bbb{F}\left( d_{1},\ldots ,d_{k}\right) $ of $A_{l}$ is
orientable if and only if $d_{1},\ldots ,d_{k},d_{k+1}$ satisfy the ${\rm %
mod}2$ condition. Here we write $d_{k+1}=n=l+1$. Alternatively orientability
holds if and only if the sizes of the blocks in $K_{\Theta }$ are congruent $%
{\rm mod}2$.
\end{proposicao}

\begin{profe}
By the comments above, the simple roots outside $\Theta $ are $\alpha
_{r_{1}},\ldots ,\alpha _{r_{k}}$, where $d_{1},\ldots ,d_{k}$ are the
dimensions determining the flag. For an index $i$ there either $%
d_{i+1}=d_{i}+1$ or $d_{i+1}>d_{i}+1$. In the second case the set $\Delta
=\{\alpha _{r_{i}+1},\ldots ,\alpha _{r_{i+1}-1}\}$ is a connected component
of $\Theta $, having $d_{i+1}-d_{i}-1$ elements. We consider two cases:

\begin{enumerate}
\item  If the second case holds for every $\alpha \notin \Theta $ then the
connected components of $\Sigma \setminus \Theta $ are singletons. If this
holds and $\alpha \notin \Theta $ is not one of the extreme roots $\alpha
_{1}$ or $\alpha _{l}$ then $\alpha $ is linked to exactly two connected
components of $\Theta $. By the first row of table \ref{tabaele} these
connected components of $\Theta $ must have the same ${\rm mod}2$ number of
elements if $\Bbb{F}\left( d_{1},\ldots ,d_{k}\right) $ is to be orientable.
Hence if $\{\alpha _{1},\alpha _{l}\}\subset \Theta $ then $\Bbb{F}\left(
d_{1},\ldots ,d_{k}\right) $ is orientable if and only if the number of
elements in the components of $\Theta $ are ${\rm mod}2$ congruent. This is
the same as the condition in the statement because a connected component has 
$d_{i+1}-d_{i}-1$ elements. On the other hand if $\alpha _{1}$ or $\alpha
_{l}$ is not in $\Theta $ then orientability holds if and only if all the
number of elements of the components of $\Theta $ are even. In this case $%
d_{i+1}-d_{i}$ is odd and $d_{1}-d_{0}=1$ or $d_{k+1}-d_{k}=1$. Hence the
result follows.

\item  As in the first case one can see that if some of the components of $%
\Sigma \setminus \Theta $ is not a singleton then all the components of $%
\Theta $ must have an even number of elements. Therefore the integers $%
d_{i+1}-d_{i}$ are odd.
\end{enumerate}
\end{profe}

\vspace{12pt}%

\noindent \textbf{Example:} A Grassmannian $\mathrm{Gr}_{k}\left( n\right) $
of $k$-dimensional subspaces in $\Bbb{R}^{n}$ is orientable if and only if $%
n $ is even.

\vspace{12pt}%

\noindent \textbf{Remark:} The orientability of the flag manifolds of 
\textrm{Sl}$\left( n,\Bbb{R}\right) $ can be decide also via Stiefel-Whitney
classes as in Conde \cite{cond}.

\subsubsection{$B_{l}=\frak{so}\left( l,l+1\right) $}

Here the flag manifolds are $\Bbb{F}^{I}\left( d_{1},\ldots ,d_{k}\right) =%
\Bbb{F}_{\Theta }$ such that $j\in \{d_{1},\ldots ,d_{k}\}$ if and only if $%
j $ is the index of a simple root $\alpha _{j}\notin \Theta $. The subgroup $%
K_{\Theta }$ is a product $\mathrm{SO}\left( n_{1}\right) \times \cdots
\times \mathrm{SO}\left( n_{s}\right) $ with the sizes $n_{i}$ given as
follows:

\begin{enumerate}
\item  If $d_{k}=l$, or equivalently $\alpha _{l}\notin \Theta $ then $%
K_{\Theta }=\mathrm{SO}\left( d_{1}\right) \times \cdots \times \mathrm{SO}%
\left( d_{k-1}-d_{k-2}\right) $.

\item  If $d_{k}<l$, or equivalently $\alpha _{l}\in \Theta $ then

\begin{enumerate}
\item  $K_{\Theta }=\mathrm{SO}\left( d_{1}\right) \times \cdots \times 
\mathrm{SO}\left( d_{k}-d_{k-1}\right) \times \mathrm{SO}\left( 2\right) $
if $d_{k}=l-1$, that is, $\{\alpha _{l}\}$ is a connected component of $%
\Theta $.

\item  $K_{\Theta }=\mathrm{SO}\left( d_{1}\right) \times \cdots \times 
\mathrm{SO}\left( d_{k}-d_{k-1}\right) \times \mathrm{SO}\left(
l-d_{k}\right) \times \mathrm{SO}\left( l-d_{k}+1\right) $ if $d_{k}<l-1$,
that is, the connected component of $\Theta $ containing $\alpha _{l}$ is a $%
B_{l-d_{k}}$.
\end{enumerate}
\end{enumerate}

\begin{proposicao}
The following two cases give necessary and sufficient conditions for flag
manifold $\Bbb{F}^{I}\left( d_{1},\ldots ,d_{k}\right) $ of $B_{l}$ to be
orientable.

\begin{enumerate}
\item  Suppose that $d_{k}=l$, that is, $\alpha _{l}\notin \Theta $. Then $%
\Bbb{F}^{I}\left( d_{1},\ldots ,d_{k}\right) $ is orientable if and only if $%
d_{1},\ldots ,d_{k-1}$, up to $k-1$, satisfy the ${\rm mod}2$ condition.
Equivalently, the sizes of the $\mathrm{SO}\left( n_{i}\right) $-components
of $K_{\Theta }$ are congruent ${\rm mod}2$.

\item  Suppose that $d_{k}<l$, that is, $\alpha _{l}\in \Theta $. Then $\Bbb{%
F}^{I}\left( d_{1},\ldots ,d_{k}\right) $ is orientable if and only if $%
d_{1},\ldots ,d_{k}$ together with $l-d_{k}$ satisfy the ${\rm mod}2$
condition.
\end{enumerate}
\end{proposicao}

\begin{profe}
If $\alpha _{l}\notin \Theta $ then $\Theta $ is contained in the $A_{l-1}$%
-subdiagram $\{\alpha _{1},\ldots ,\alpha _{l-1}\}$. Hence the condition is
the same as in the $A_{l}$ case. Furthermore, $S\left( \alpha _{l},\Delta
\right) $ is even for any $\Delta $ because $\alpha _{l}$ is a short root.
Therefore no further condition comes in.

In the second case, if $\Delta $ is the connected component of $\Theta $
containing $\alpha _{l}$ then the contribution $S\left( \alpha ,\Delta
\right) $ of $\Delta $ to the total sum is the number of elements of $\Delta 
$ by tables \ref{tabaele} and \ref{tabbele}. Again, the conclusion is as in
the $A_{l}$ case.
\end{profe}

\vspace{12pt}%

\noindent \textbf{Example:} A Grassmannian $\mathrm{Gr}_{k}^{I}\left(
n\right) =\Bbb{F}^{I}\left( k\right) $ of $k$-dimensional isotropic
subspaces in $\Bbb{R}^{2l+1}$ is orientable if and only if either i) $k=l$
or ii) $k<l$ and $l$ is even.

\subsubsection{$C_{l}=\frak{sp}\left( l,\Bbb{R}\right) $}

Again the flag manifolds are $\Bbb{F}^{I}\left( d_{1},\ldots ,d_{k}\right) =%
\Bbb{F}_{\Theta }$ such that $j\in \{d_{1},\ldots ,d_{k}\}$ if and only if $%
j $ is the index of a simple root $\alpha _{j}\notin \Theta $. The subgroup $%
K_{\Theta }$ is

\begin{enumerate}
\item  $\mathrm{SO}\left( d_{1}\right) \times \cdots \times \mathrm{SO}%
\left( d_{k-1}-d_{k-2}\right) $ if $d_{k}=l$.

\item  $\mathrm{SO}\left( d_{1}\right) \times \cdots \times \mathrm{SO}%
\left( d_{k-1}-d_{k-2}\right) \times \mathrm{SO}\left( 2\right) $ if $%
d_{k}=l-1$.

\item  $\mathrm{SO}\left( d_{1}\right) \times \cdots \times \mathrm{SO}%
\left( d_{k-1}-d_{k-2}\right) \times \mathrm{U}\left( l-d_{k}\right) $ if $%
d_{k}<l-1$.
\end{enumerate}

\begin{proposicao}
For $C_{l}$ a necessary and sufficitent condition for the orientability of $%
\Bbb{F}^{I}\left( d_{1},\ldots ,d_{k}\right) $ is that $d_{1},\ldots ,d_{k}$
satisfy the ${\rm mod}2$ condition.
\end{proposicao}

\begin{profe}
There are two possibilities:

\begin{enumerate}
\item  If $d_{k}=l$, that is, $\alpha _{l}\notin \Theta $ then $\Theta $ is
contained in the $A_{l-1}$ and the condition, up to $k-2$, comes from the $%
A_{l}$ case. The difference $d_{k}-d_{k-1}$ also enters in the condition
because $\alpha _{l}$ is a large root.

\item  If $d_{k}<l$, that is, $\alpha _{l}\in \Theta $ then the conditions
are necessary as in the $A_{l}$ case. To see that no further condition
appears look at the connected component $\Delta $ containing $\alpha _{l}$.
If $\Delta =\{\alpha _{l}\}$ then $S\left( \alpha _{l-1},\Delta \right) $ is
even because $\alpha _{l-1}$ is a short root. Otherwise, $\Delta $ is a $%
C_{k}$ and its contribution is also even by table \ref{tabcele}.
\end{enumerate}
\end{profe}

\subsubsection{$D_{l}=\frak{so}\left( l,l\right) $}

The flag manifolds of $\frak{so}\left( l,l\right) $ are also realized as
flags of isotropic subspaces with a slight difference from the odd
dimensional case $B_{l}=\mathrm{SO}\left( l,l+1\right) $. First a minimal
flag manifold $\Bbb{F}_{\Sigma \setminus \{\alpha _{i}\}}$ is the
Grassmannian of isotropic subspaces of dimension $i$ if $i\leq l-2$.
However, both $\Bbb{F}_{\Sigma \setminus \{\alpha _{l-1}\}}$ and $\Bbb{F}%
_{\Sigma \setminus \{\alpha _{l}\}}$ are realized as subsets of $l$%
-dimensional isotropic subspaces. Each one is a closed orbit of the identity
component of $\mathrm{SO}\left( l,l\right) $ in the Grassmannian $\mathrm{Gr}%
_{l}^{I}\left( 2l\right) $ of $l$-dimensional isotropic subspaces. We denote
these orbits by $\mathrm{Gr}_{l^{+}}^{I}\left( 2l\right) =\Bbb{F}_{\Sigma
\setminus \{\alpha _{l}\}}$ and $\mathrm{Gr}_{l^{-}}^{I}\left( 2l\right) =%
\Bbb{F}_{\Sigma \setminus \{\alpha _{l}-1\}}$. (By the way the isotropic
Grassmannian $\mathrm{Gr}_{l-1}^{I}\left( 2l\right) $ is the flag manifold $%
\Bbb{F}_{\Sigma \setminus \{\alpha _{l-1},\alpha _{l}\}}$, which is not
minimal.)

Accordingly, the flag manifolds of $\frak{so}\left( l,l\right) $ are defined
by indices $1\leq d_{1}\leq \cdots \leq d_{k}\leq l-2$ joined eventually to $%
l^{+}$ and $l^{-}$. The elements of $\Bbb{F}^{I}\left( d_{1}\ldots
,d_{k}\right) $ are flags of isotropic subspaces $V_{1}\subset \cdots
\subset V_{k}$ with $\dim V_{i}=d_{k}$. When $l^{+}$ or $l^{-}$ are present
then one must include an isotropic subspace in $\mathrm{Gr}%
_{l^{+}}^{I}\left( 2l\right) $ or $\mathrm{Gr}_{l^{-}}^{I}\left( 2l\right) $%
, respectively, containing $V_{k}$, and hence the other subspaces.

The group $K_{\Theta }$ is a product of $\mathrm{SO}\left( d\right) $'s
components each one for a connected component of $\Theta $ unless a $D_{k}$
component appears. Such a component contributes to $K_{\Theta }$ with a $%
\mathrm{SO}\left( k\right) \times \mathrm{SO}\left( k\right) $.

\begin{proposicao}
The orientability of the flag manifolds of $D_{l}=\frak{so}\left( l,l\right) 
$ is given as follows:

\begin{enumerate}
\item  For a flag $\Bbb{F}^{I}\left( d_{1},\ldots ,d_{k}\right) $ there are
the possibilities:

\begin{enumerate}
\item  If $d_{k}\leq l-4$ then orientability holds if and only if $%
d_{1},\ldots ,d_{k}$ satisfy the ${\rm mod}2$ condition.

\item  If $d_{k}=l-3$ then orientability holds if and only if the
differences $d_{i+1}-d_{i}$, $i=0,\ldots ,k-1$, are even numbers.

\item  If $d_{k}=l-2$ then orientability holds if and only if the
differences $d_{i+1}-d_{i}$, $i=0,\ldots ,k-1$, are odd numbers.
\end{enumerate}

\item  For the flag manifolds $\Bbb{F}^{I}\left( d_{1},\ldots
,d_{k},l^{+}\right) $ and $\Bbb{F}^{I}\left( d_{1},\ldots
,d_{k},l^{-}\right) $ we have:

\begin{enumerate}
\item  If $d_{k}=l-2$ then the condition is that $d_{i+1}-d_{i}$, $%
i=0,\ldots ,k-2$, are even numbers.

\item  If $d_{k}<l-2$ then the condition is that $d_{i+1}-d_{i}$, $%
i=0,\ldots ,k-2$, are odd numbers and $d_{k}-d_{k-1}$ is even.
\end{enumerate}

\item  For the flag manifolds $\Bbb{F}^{I}\left( d_{1},\ldots
,d_{k},l^{+},l^{-}\right) $ we have:

\begin{enumerate}
\item  If $d_{k}=l-2$ then $d_{1},\ldots ,d_{k-2}$ satisfy the ${\rm mod}2$
condition.

\item  If $d_{k}<l-2$ then $d_{i+1}-d_{i}$, $i=0,\ldots ,k-2$, are odd
numbers.
\end{enumerate}
\end{enumerate}
\end{proposicao}

\begin{profe}
If $d_{k}\leq l-4$ then $\Theta $ contains a connected component $\Delta $
which is a $D_{k}$ (at the right side of the diagram). By table \ref{tabdele}
the contribution of $\Delta $ is even, so that orientability depends on the
roots in the $A_{l-4}$ diagram $\{\alpha _{1},\ldots ,\alpha _{l-4}\}$ where
the condition is as in the statement. If $d_{k}=l-3$ then the differences $%
d_{i+1}-d_{i}$, $i=0,\ldots ,k-1$, must be congruent ${\rm mod}2$ to have
orientability. But the root $\alpha _{l-3}$ is linked to the $A_{3}=\{\alpha
_{l-2},\alpha _{l-1},\alpha _{l}\}$, so that the number of elements of the
components of $\Theta $ are odd, that is, the differences $d_{i+1}-d_{i}$
are even. The same argument applies to $d_{k}=l-2$ , but now $\alpha _{l-2}$
is linked to the two $A_{1}$'s $\{\alpha _{l-1}\}$ and $\{\alpha _{l}\}$.

The other cases are checked the same way.
\end{profe}

\section{ Vector bundles over flag bundles}

\label{secbundlesbundles}

In this final section we consider vector bundles over flag bundles. The
orientability of vector bundles over the flag manifolds carry over to vector
bundles over flag bundles in case the latter are bundles associated to
trivial principal bundles.

With the previous notation let $R$ a $K$-principal bundle. Since $K$ acts
continuously on $V$ and $X$, the associated bundle $R\times _{K}V$ is a
finite dimensional vector bundle over $R\times _{K}X$ whose fibers are the
same as the fibers of $V$.

\begin{proposicao}
\label{proporientassociado}Assume that $R$ is trivial. Then the vector
bundle 
\[
R\times _{K}V\to R\times _{K}X
\]
is orientable if, and only if, the vector bundle $V\to X$ is orientable.
\end{proposicao}

\begin{profe}
Since the $K$-principal bundle $R\to Y$ is trivial, we have that $R\times
_{K}V\to R\times _{K}X$ is homeomorphic as a vector bundle to $Y\times V\to
Y\times X$. Since the frame bundle of $Y\times V$ can be given by $Y\times BV
$, the orientation bundle of $Y\times V$ can be given by $Y\times \mathcal{O}%
V$. If $\sigma :X\to \mathcal{O}V$ is a continuous section, then $%
(y,x)\mapsto (y,\sigma (x))$ is a continuous section of $Y\times \mathcal{O}V
$. Reciprocally, if $\sigma :Y\times X\to Y\times \mathcal{O}V$ is a
continuous section, then $x\mapsto \sigma (y_{0},x)$ is a continuous section
of $\mathcal{O}V$, where $y_{0}\in Y$.
\end{profe}

Let $G$ be a Lie group acting on its Lie algebra $\mathfrak{g}$ by the
adjoint action. The vector bundles we will consider in the sequel arise as
associated bundles of the $L$-principal bundle $K \to K/L$, where $K$ is a
subgroup of $G$. For an $L$-invariant subspace $\mathfrak{l}$ of $%
\mathfrak{g}$, we will consider the associated vector bundle 
\[
V = K \times_L \mathfrak{l}, 
\]
whose typical fiber is $\mathfrak{l}$.

\begin{corolario}
\label{cororient}The associated vector bundle $V$ is orientable if and only
if $\det (g|_{\mathfrak{l}})>0$, for every $g\in L$.
\end{corolario}

\begin{profe}
We only need to show that $V$ satisfies the hypothesis of Proposition \ref
{proporient}. First we note that its frame bundle is given by $BV=K\times
_{L}\mathrm{Gl}(\mathfrak{l})$. Defining an action $k\in K$ on $m\cdot X\in
BV$ by 
\[
k(m\cdot X)=km\cdot X,
\]
where $m\in K$, $X\in \mathfrak{l}$, we have that the action of $K$ on $K/L$
lifts to a continuous action of automorphisms on the frame bundle $BV$.
\end{profe}

To conclude we apply our results to the situation of \cite{psmsconley},
where flows on flag bundles and their Conley indices are considered. In \cite
{psmsconley} one starts with a principal bundle $Q\rightarrow X$ whose
structural group $G$ is semi-simple, and a flow $\phi _{t}$, $t\in \Bbb{Z}$
or $\Bbb{R}$, of automorphisms of $Q$. There are induced flows on the
associated bundles $Q\times _{G}F$, where the typical fiber $F$ is acted by $%
G$ on the left. In particular, in \cite{psmsconley} it is taken as a typical
fiber $F$ a flag manifold $\Bbb{F}_{\Theta }$ of $G$ yielding the flag
bundle $\Bbb{E}_{\Theta }=Q\times _{G}\Bbb{F}_{\Theta }$.

According to the results of \cite{msm} and \cite{psmsconley}, each Morse
component $\mathcal{M}_{\Theta }(w)$ of $\phi ^{t}$ is a flag bundle of a
certain subbundle $Q_{\phi }$ of $Q$. Moreover, the unstable set $\mathcal{V}%
_{\Theta }^{+}(w)$ of the Morse component $\mathcal{M}_{\Theta }(w)$ is an
associated vector bundle of $Q_{\phi }$ whose base is $\mathcal{M}_{\Theta
}(w)$ and whose typical fiber is the same as the fiber of $V_{\Theta
}^{+}(H_{\phi },w)$, where is $H_{\phi }$ is a certain element of $\mathrm{cl%
}\frak{a}^{+}$, called the parabolic type of $\phi ^{t}$.

When the base $B$ is a point, the flow of automorphisms $\phi ^{t}$ is given
by $g^{t}$ for some $g\in G$, when $t\in \Bbb{Z}$, or by $\exp (tX)$ for
some $X\in \mathfrak{g}$, when $t\in \Bbb{R}$. In \cite{fps}, it is shown
that the parabolic type $H_{\phi }$ of these flows is given by the
hyperbolic component of $g$ or $X$ under the Jordan decomposition.

In \cite{psmsconley}, we show that the Conley index of the attractor
component in the maximal flag bundle and, under certain hypothesis, the
Conley index of each Morse component, is the Thom space of its unstable
vector bundle. The orientability of the unstable vector bundle then comes to
the scene in order to apply Thom isomorphism and detect the homological
Conley indices of the Morse components. With these results in mind we state
the following criterion of orientability of $\mathcal{V}_{\Theta }^{+}(w)$,
that follows immediately from Proposition \ref{proporientassociado}.

\begin{proposicao}
\label{propfibradoinstavel}Assume that the reduction $R_{\phi }$ is a
trivial bundle. The stable and unstable vector bundles $\mathcal{V}_{\Theta
}^{\pm }(H,w)$ are orientable if and only if the vector bundles $V_{\Theta
}^{\pm }(H,w)$ are orientable.
\end{proposicao}

There are two cases where  the hypothesis of the above result are
automatically satisfied. Namely for periodic flows, it is shown in \cite{fps}
that the reduction $Q_{\phi }$ is trivial. For the control flow of \cite{ck}, the reduction $Q_{\phi }$ is always trivial since the base space of the
control flow is contractible.

\end{document}